\documentclass[10pt]{article}

\usepackage{a4wide}
\usepackage{amssymb}
\usepackage{amsfonts}
\usepackage{amsmath}
\input xy
\xyoption{arrow}
\xyoption{matrix}

\date{}

\newtheorem{proposition}{Proposition}[section]
\newtheorem{theorem}[proposition]{Theorem}
\newtheorem{lemma}[proposition]{Lemma}

\newtheorem{corollary}[proposition]{Corollary}

\def\GK{{\rm  GK}\,}
\def\Kdim{{\rm K.dim }\,}

\def\der{\partial }

\def\nFM0{{\nu }_{F,M_0}}
\def\nFN0{{\nu }_{F,N_0}}
\def\nGN0{{\nu }_{G,N_0}}

\def\N0{ {\bf N}_0 }

\def\ra{\rightarrow}

\def\lra{\leftrightarrow}
\def\Xpm{X^{\pm }}

\def\s{\sigma}
\def\Z{{\mathbb Z }}

\def\l1{{\lambda}_1}

\def\a{\alpha}
\def\a0{ {\alpha }_0}
\def\a1{ {\alpha }_1}

\def\l{\lambda}
\def\o{\omega}

\def\nFGM0{{\nu }_{F,G,M_0}}

\def\nFN0{{\nu}_{F,N_0}}

\def\sm{{\sigma}^m}

\def\sm1{{\sigma}^{-1}}

\def\smtp1{{\sigma}^{-t+1}}

\def\o{\omega }
\def\S1{S^{-1}}

\def\Xpm1{X^{\pm 1}_1}

\def\sPM1{{\sigma }^{\pm 1}}
\def\sMP1{{\sigma }^{\mp 1 }}

\def\d{\delta}

\def\L{\Lambda}

\def\OO{{\cal O}}
\def\CA{{\cal A}}

\def\CD{{\cal D}}

\def\Ytm1{Y^{t-1}}
\def\Yim1{Y^{i-1}}

\def\i{{\bf i}}

\def\Der{{\rm Der }}

\def\D{ \Delta }

\def\SL2Z{ {\rm SL}_2({\bf Z}) }

\def\Gp1{ G^{1 , 1 } }
\def\P11{ P^{-1 , 1 } }
\def\Pp1{ P^{1 , 1 } }

\def\nCLsr{{}^\nu\kern-2pt {\cal L}^{\sigma , \rho  }}
\def\nP{{}^\nu \kern-2pt P}
\def\nL{{}^\nu\kern-2pt L}
\def\nLL{{}^\nu\kern-2pt \Lambda}
\def\nPsr{{}^\nu\kern-2pt P^{\sigma , \rho  }}
\def\nLsr{{}^\nu\kern-2pt L^{\sigma , \rho  }}
\def\nuCL{{}^\nu\kern-2pt  {\cal L}}
\def\nCLsr{{}^\nu\kern-2pt {\cal L}^{\sigma , \rho  }}
\def\nCL1m{{}^\nu\kern-2pt {\cal L}^{-1 , 1  }}
\def\x1nu{x^\frac{1}{\nu}}
\def\xm1nu{x^{-\frac{1}{\nu}}}

\def\ra{\rightarrow }

\def\CC{ {\cal C}}

\def\nAM0{{\nu }_{{\cal A},M_0}}
\def\nAN0{{\nu }_{{\cal A},N_0}}

\def\Kdim{ {\rm Kdim } }

\def\End{ {\rm End }}
\def\Der{ {\rm Der }}

\def\det{ {\rm det }}

\def\bp{\overline{p}}

\def\ga{\mathfrak{a}}

\def\bJ{\overline{J}}

\def\j{{\bf j}}

\def\II{{\bf I}}
\def\JJ{{\bf J}}

\begin{document}

\author{V. V. \  Bavula   
}

\title{Simplicity criteria for  rings of differential operators}

\maketitle
\begin{abstract}

Let $K$ be a field of arbitrary characteristic, $\CA$ be a commutative $K$-algebra which is a domain of essentially finite type (eg, the algebra of functions on an irreducible affine algebraic variety), $\ga_r$ be its {\em Jacobian ideal}, $\CD (\CA )$ be the algebra of differential operators on the algebra $\CA$. The aim of the paper is to give a simplicity criterion for the algebra $\CD (\CA )$: {\em The algebra $\CD (\CA )$ is simple iff $\CD (\CA ) \ga_r^i\CD (\CA )= \CD (\CA )$ for all $i\geq 1$ provided the field $K$  is a perfect field.} Furthermore, a simplicity criterion is given for the algebra $\CD (R)$ of differential operators on an arbitrary commutative algebra $R$ over an arbitrary field.  This gives an answer to an old question to find a simplicity criterion for algebras of differential operators. \\ 

\end{abstract}


\section{Introduction}

The following notation will remain  fixed throughout the paper (if
it is not stated otherwise): $K$ is a field of arbitrary  characteristic 
(not necessarily algebraically closed), module means a left
module, $P_n=K[x_1, \ldots , x_n]$ is a polynomial algebra over
$K$, $\der_1:=\frac{\der}{\der x_1}, \ldots ,
\der_n:=\frac{\der}{\der x_n}\in \Der_K(P_n)$, $I:=\sum_{i=1}^m
P_nf_i$ is a {\bf prime} but {\bf not} a maximal ideal of the
polynomial algebra $P_n$ with a set of generators $f_1, \ldots ,
f_m$, the algebra $A:=P_n/I$ which is a domain with the field of
fractions $Q:={\rm Frac}(A)$, the epimorphism $\pi :P_n\ra A$,
$p\mapsto \bp :=p+I$, to make notation simpler we sometime write
$x_i$ for $\overline{x}_i$ (if it does not lead to confusion), the
{\bf Jacobi} $m\times n$ matrices
 $J=(\frac{\der f_i}{\der x_j})\in M_{m,n}(P_n)$ and
 $\bJ =(\overline{\frac{\der f_i}{\der
x_j}})\in M_{m,n}(A)\subseteq M_{m,n}(Q)$,  $r:={\rm rk}_Q(\bJ )$
is the {\bf rank} of the Jacobi matrix $\bJ$ over the field $Q$,
$\ga_r$ is the {\bf Jacobian ideal} of the algebra $A$ which is
(by definition) generated by all the $r\times r$ minors of the
Jacobi matrix $\bJ$ (Suppose that $K$ is a perfect field. Then the algebra $A$ is {\em regular} iff $\ga_r=A$, it is the
{\bf Jacobian criterion of regularity}, \cite[Theorem 16.19]{Eisenbook}).
 For $\i =(i_1, \ldots , i_r)$ such that $1\leq
i_1<\cdots <i_r\leq m$ and $\j =(j_1, \ldots , j_r)$ such that
$1\leq j_1<\cdots <j_r\leq n$, $\D
 (\i , \j )$ denotes the corresponding minor of the Jacobi matrix
$\bJ =(\bJ_{ij})$, that is $\det (\bJ_{i_\nu, j_\mu})$, $\nu , \mu
=1, \ldots, r$,  and the element $\i$ (resp., $\j $) is called {\bf
non-singular} if $\D (\i , \j')\neq 0$ (resp., $\D (\i', \j )\neq
0$) for some $\j'$ (resp., $\i'$). We denote by $\II_r$ (resp.,
$\JJ_r$) the set of all the non-singular $r$-tuples $\i$ (resp.,
$\j $).

Since $r$ is the rank of the Jacobi matrix $\bJ$, it is easy to
show that $\D (\i , \j )\neq 0$ iff $\i\in \II_r$ and $\j\in
\JJ_r$, \cite[Lemma 2.1]{gendifreg}. 

 A localization of an {\em affine} algebra is
called an algebra of {\bf essentially finite type}. Let $\CA :=\S1 A$ be a localization of the algebra $A=P_n/I$ at a
multiplicatively closed subset $S$ of $A$. Suppose that $K$ is a perfect field. Then the algebra $\CA$ is {\em regular} iff $\ga_r=\CA$ where $\ga_r$ is the Jacobian ideal of $\CA$, it is the
{\bf Jacobian criterion of regularity}, \cite[Theorem 16.19]{Eisenbook}.
For any regular algebra $\CA$ over a perfect field,  explicit sets of generators and defining relations for the algebra $\CD (\CA )$ are given in \cite{gendifreg} (char($K$)=0) and  \cite{gendifregcharp} (char($K)>0$).

Let $R$ be an arbitrary commutative $K$-algebra. We denote by $\CD (R)$ the algebra of differential operators on the algebra $R$ and by 
$\Der_K(R)$ the $R$-module of $K$-derivations of
 $R$. The action of a derivation $\d$ on an element $a$
is denoted by $\d (a)$. \\

{\bf Simplicity criterion for the algebra $\CD (\CA )$ where the algebra $\CA$ is a domain of essentially finite type.} Theorem \ref{9Jul19} is a simplicity criterion for the algebra $\CD (\CA )$ where the algebra $\CA$ is a domain of essentially finite type.

\begin{theorem}\label{9Jul19}
Let  a $K$-algebra $\CA$ be a commutative domain of essentially finite type over a perfect field $K$, and $\ga_r$ be its Jacobial ideal. The following statements are equivalent:
\begin{enumerate}
\item The algebra $\CD (\CA )$ of differential operators on $\CA$ is a simple algebra.
\item For all $i\geq 1$, $\CD (\CA ) \ga_r^i\CD (\CA )= \CD (\CA )$. 
\item For all $k\geq 1$, $\i\in \II_r$ and $\j\in \JJ_r$,  $\CD (\CA ) \D (\i, \j )^k\CD (\CA )= \CD (\CA )$.
\end{enumerate}
\end{theorem}
As an application of Theorem \ref{9Jul19} we show that the algebra of differential operators on the cusp is simple.\\

{\bf Simplicity criterion for the algebra $\CD (R )$ where $R$ is an arbitrary commutative algebra.} An ideal $\ga $ of the algebra $R$ is called $\Der_K(R)$-{\em stable} if $\d (\ga ) \subseteq \ga$ for all $\d \in \Der_K(R)$. Theorem \ref{A9Jul19}.(2) is a simplicity criterion for the algebra $\CD (R )$ where $R$ is an arbitrary commutative algebra. Theorem \ref{A9Jul19}.(1) shows that every nonzero ideal of the algebra $\CD (R)$ meets the subalgebra $R$ of $\CD (R)$. If, in addition, the algebra $R=\CA $ is a domain of essentially finite type,  Theorem \ref{A9Jul19}.(3) shows that every nonzero ideal of the algebra $\CD (R)$ contains a power of the Jacobian ideal of $\CA$.

\begin{theorem}\label{A9Jul19}
Let $R$ be a commutative algebra over an arbitrary field $K$. 
\begin{enumerate}
\item Let $I$ be a nonzero ideal the algebra $\CD (R)$. Then the ideal $I_0:= I\cap R$ is a nonzero $\Der_K(R)$-stable ideal of the algebra $R$ such that $\CD (R)I_0\CD (R)\cap R = I_0$. In particular, every nonzero ideal of the algebra $\CD (R)$ has nonzero intersection with $R$.
\item The ring $\CD (R)$ is not simple iff there is a proper $\Der_K(R)$-stable ideal $\ga $ of $R$ such that $\CD (R)\ga \CD (R)\cap R = \ga$.
\item Suppose, in addition, that $K$ is a perfect field and  the algebra $\CA = R$ is a domain of essentially finite type, $\ga_r$ be its Jacobian ideal, $I$ be a nonzero ideal of $\CD (\CA )$ and $I_0=I\cap \CA$. Then $\ga_r^i\subseteq I_0$ for some $i\geq 1$. 
\end{enumerate}
\end{theorem}


\section{Proofs of Theorem \ref{9Jul19} and Theorem \ref{A9Jul19}}\label{PRFS}

In this section, proofs of Theorem \ref{9Jul19} and Theorem \ref{A9Jul19} are given.

Let $R$ be a commutative $K$-algebra. The ring of ($K$-linear)
{\bf differential operators} $\CD (R)$ on $R$ is defined as a
union of $R$-modules  $\CD (R)=\bigcup_{i=0}^\infty \,\CD_i (R)$
where 
$$ \CD_i (R)=\{ u\in \End_K(R)\, |\, [r,u]:=ru-ur\in \CD_{i-1} (R)\; {\rm for \; all \; }\; r\in R\},\;\; i\geq 0, \;\; \CD_{-1}(R):=0.$$
In particular, $\CD_0 (R)=\End_R(R)\simeq R$, $(x\mapsto bx)\lra b$.
 The set of $R$-bimodules $\{ \CD_i (R)\}_{i\geq 0}$ is the {\bf order filtration} for
the algebra $\CD (R)$:
$$\CD_0(R)\subseteq   \CD_1 (R)\subseteq \cdots \subseteq
\CD_i (R)\subseteq \cdots\;\; {\rm and}\;\; \CD_i (R)\CD_j
(R)\subseteq \CD_{i+j} (R) \;\; {\rm for\; all}\;\; i,j\geq 0.$$

The subalgebra $\D (R)$ of $\CD (R)$ which is  generated by $R\equiv
\End_R(R)$ and the set ${\rm Der}_K (R)$ of all $K$-derivations of
$R$ is called the {\bf derivation ring} of $R$.

Suppose that $R$ is a  regular affine  domain of Krull dimension
$n\geq 1 $ and char($K$)=0. In geometric terms, $R$ is the coordinate ring $\OO
(X)$ of a smooth irreducible  affine algebraic variety $X$ of
dimension $n$. Then
\begin{itemize}
\item ${\rm Der}_K (R)$ {\em is a finitely generated projective}
$R$-{\em module of rank} $n$, \item  $\CD (R)=\Delta (R) $, \item
$\CD (R)$ {\em is a simple (left and right) Noetherian domain of
Gelfand-Kirillov dimension}  $\GK \, \CD (R)=2n$ ($n=\GK (R)=\Kdim
(R))$.
\end{itemize}

For the proofs of the statements above the reader is referred to
\cite{MR}, Chapter 15.
 So, the domain $\CD (R)$ is a simple finitely generated infinite dimensional Noetherian algebra
(\cite{MR}, Chapter 15).  

If char($K)>0$ then  $\CD (R)\neq\Delta (R) $ and the algebra $\CD (R$ is not finitely generated and neither left nor right Noetherian but analogues of the results above hold but the Gelfand-Kirillov dimension has to replaced by a new dimension introduced in \cite{BernIneqcharp}.

Given a ring $B$ and a non-nilpotent element $s\in B$. Suppose that the set $S_s:=\{ s^i\, | \, i\geq 0\}$ is a left denominator set of $B$. The localization $S_s^{-1}B$ of the ring $B$ at $S_s$ is also denoted by $B_s$.\\

{\bf Proof of Theorem \ref{A9Jul19}}. 1. (i) {\em The ideal $I_0$ of $R$ is a $\Der_K(R)$-stable ideal}: For all $\d \in \Der_K(R)$, $I_0\supseteq [\d , I_0]=\d (I_0)$.

(ii) $\CD (R)I_0\CD (R)\cap R = I_0$: $I_0 \subseteq \CD (R)I_0\CD (R)\cap R \subseteq  \CD (R)I\CD (R)\cap R = I\cap R=I_0$, and the statement  (ii) follows. 

(iii) $I_0\neq 0$: Recall that the ring $\CD (R)$ admits the order filtration $\{ \CD (R)_i\}_{i\geq 0}$. Therefore, $I=\bigcup_{i\geq 0}I_i$ where $I_i= I\cap \CD (R)_i$. Let $s=\min \{ i\geq 0\, | \, I_i\neq 0\}$. Then $I_s\neq 0$ and 
$$[r,I_s]\subseteq I_{s-1}=\{ 0\}=\CD_{-1}(R)\;\; {\rm for\; all} \;\; r\in R, $$
i.e. $I_s \subseteq \CD (R)_0=R$, by the {\em definition of the order filtration} on $\CD (R)$, and so $s=0$, as required. 

2. $(\Rightarrow )$ If the ring $\CD (R)$ is not simple then there is  proper ideal, say $I$, of $\CD (R)$. Then, by the statements (i) and (ii) in the proof of statement 1, it suffices to take $\ga = I_0$.

$(\Leftarrow )$ The implication is obvious. 

3. Recall that the Jacobian ideal $\ga_r$ of the algebra $\CA$ is generated by the {\em finite} set $\{ \D (\i , \j )\, | \, \i\in \II_r, \j\in \JJ_r\}$. For each element $\D (\i , \j )$, the algebra $\CA_{\D (\i , \j )}$ is a regular domain of essentially finite type. So, the algebra $\CD (\CA_{\D (\i , \j )})\simeq \CD (\CA)_{\D (\i , \j )}$ is simple  (the algebra $\CD (\CA)_{\D (\i , \j )}$ is a left and right localization of $\CD (\CA)$ at the powers of the element $\D (\i , \j )$). Therefore, $1\in I_{\D (\i , \j )}$, and so $\D (\i , \j )^l\in I\cap \CA =I_0$ for some $l\geq 1$. So, $\ga_r^i\subseteq I_0$ for some $i\geq 1$.  $\Box $\\

{\bf Proof of Theorem \ref{9Jul19}}. $(1\Rightarrow 3)$ The implication is trivial.
 
$(3\Rightarrow 2)$ The implication follows from the fact that the Jacobian ideal $\ga_r$ of the algebra $\CA$ is generated by the finite set $\{ \D (\i , \j )\, | \, \i \in \II_r, \j\in \JJ_r\}$. In particular, $\D (\i , \j )^k\subseteq \ga_r^k$ for all $k\geq 1$, and so $\CD (\CA )=\CD (\CA ) \D (\i, \j )^k\CD (\CA )\subseteq \CD (\CA ) \ga_r^k\CD (\CA )$.

$(2\Rightarrow 1)$ Suppose that the algebra $\CD (\CA )$ is not simple, we seek a contradiction. Fix a proper ideal, say $I$, of the algebra $\CD (\CA )$. By Theorem \ref{A9Jul19}.(3), $\ga_r^i\subseteq I_0$ for some natural number $i\geq 1$. Then 
$$ \CA \neq I_0=\CD (\CA ) I_0\CD (\CA)\cap \CA \supseteq  \CD (\CA ) \ga_r^i\CD (\CA)\cap \CA .$$
Therefore, $\CD (\CA ) \ga_r^i\CD (\CA)\neq \CD (\CA )$, a contradiction. $\Box $

$\noindent $

Given a commutative algebra $R$, we denote by $\CC_R$ the set of {\em regular elements} of $R$ (i.e. non-zero-divisors) and by $Q(R):=\CC_R^{-1}R$ its {\em quotient algebra}. 

\begin{corollary}\label{a11Jul19}
Let $\CA $ be a semiprime commutative algebra with finitely many minimal primes. Then the algebra $\CD (\CA )$ is a simple algebra iff the algebra $\CA$ is a domain and the algebra $\CD (\CA )$ is a simple algebra.
\end{corollary}

{\it Proof}. $(\Rightarrow )$ Suppose that the algebra $\CA$ is not a domain. Then its  quotient algebra $Q(\CA ) := \CC_\CA^{-1}\CA \simeq \prod_{i=1}^sK_i$ is a direct product of fields $K_i$ where $s\geq 2$ is the number of minimal primes of the algebra $\CA$. Therefore,
$$ \CC_\CA^{-1}(\CD (\CA ) )\simeq \CD (\CC_\CA^{-1}\CA )\simeq \CD (Q(\CA ))\simeq \CD(\prod_{i=1}^sK_i)\simeq \prod_{i=1}^s\CD (K_i).$$
The algebra $\CD (\CA )$ is an essential left $\CD (\CA )$-submodule of $\CD (Q(\CA ))$. Therefore, the intersection $ \CD (\CA )\cap \CD (K_1)$ is a proper ideal of the algebra 
 $\CD (\CA )$ since $s\geq 2$, a contradiction.
 
 $(\Leftarrow )$ The implication is trivial. $\Box $ \\


%

{\em Example.} {\sc (The algebra of differential operators on the cusp)} 
 Let $A=K[x,y]/(y^2-x^3)$, the algebra of regular functions on the cusp $y^2=x^3$. The algebra $A$ is isomorphic to the subalgebra $K+\sum_{i\geq 2}Kx^i$ of the polynomial algebra $K[x]$. Notice that $A\subseteq K[x]\subseteq A_x=K[x]_x=K[x,x^{-1}]$ and $\CD (
K[x,x^{-1}])=\oplus_{i\in \Z}Dx^i=D[x,x^{-1}; \s ]$ is a skew Laurent polynomial ring with coefficients in the polynomial algebra $D=K[h]$,  where  $h=x\der$,  and $\s$ is a $K$-automorphism of $D$ given by the rule $\s (h)=h-1$.  The algebra $A_1=K\langle x, \der \, | \, \der x -x \der =1 \rangle$ is called the (first) {\em Weyl algebra}. Then $A_1\simeq \CD (K[x])$ and $A_{1,x}\simeq \CD (K[x])_x\simeq \CD (K[x]_x)\simeq D[x,x^{-1}; \s ]\simeq \CD (A)_x$. Notice that $\CD (A) = \{ \d \in \CD(A)_x \, | \, \d (A)\subseteq A\}$.

\begin{lemma}\label{a13Jul19}
Let $A= K+\sum_{i\geq 2}Kx^i (\simeq K[x,y]/(y^2-x^3))$ and $D=K[h]$. Then 
\begin{enumerate}
\item $\CD (A)=\bigoplus_{i\in \Z}Dw_i\subseteq \CD (A_x)$ where $w_0=1$, $w_1=(h-1)x$ and $w_i=x^i$ for $i\geq 2$; $w_{-1}=(h+1)(h-1)x^{-1}$, $w_{-2} =(h+2)(h-1)x^{-2}$ and $w_{-i}=(h-1) \cdot (h+1)\cdots (h+i-2)\cdot (h+i)x^{-i}$ for $i\geq 3$.
\item The algebra $\CD (A)$ is a simple finitely generated Noetherian domain. Furthermore, the elements $h$ and  $w_i$ ($i=\pm 1, \pm 2, \pm 3$) are algebra generators of $\CD (A)$. For all $i\geq 1$, $w_{-2i}=w_{-2}^i$ and  $w_{-3-2i}=w_{-3}w_{-2}^i$.  
\item $\Der_K(A)=K[x]h$, $\D (A) = K[h][x; \s ]$ is a non-simple Noetherian algebra and $\D (A) \neq \CD (A)$. 
\item The Jacobian ideal $\ga_1=\sum_{i\geq 2}Kx^i$ of $A$ is $\D (R)$-stable but not $\CD (R)$-stable.
\end{enumerate}
\end{lemma}

{\it Proof}. 1. Recall  that $\CD (A_x)=\bigoplus_{i\in \Z}Dx^i$ is a $\Z$-graded algebra ($Dx^iDx^j\subseteq Dx^{i+j}$ for all $i,j\in \Z$), the $\CD (A_x)$-module $A_x=K[x,x^{-1}]$ is a $\Z$-graded module and the algebra $A$ is a homogeneous subalgebra of $A_x$. Now, statement 1 follows from obvious computations and the fact that $\CD (A) = \{ \d \in \CD(A)_x \, | \, \d (A)\subseteq A\}$.

2. The Jacobian ideal $\ga_1$ of $A$ is equal to $\sum_{i\geq 2}Kx^i$. Since $x^{2i}\in \ga_1^i$ for all $i\geq 1$, in order to prove simplicity of $\CD (A)$ it suffices to show that $(x^i)=\CD (A)$ for all $i\geq 2$, by Theorem \ref{9Jul19}. Notice that the polynomials of $D=K[h]$, $w_{-i}x^i$ and $x^iw_{-i}$, are coprime, hence $(x^i)=\CD (A)$ for all $i\geq 2$. In more detail, $w_{-2}x^2=(h+2)(h-1)$ and $x^2w_{-2}=h(h-3)$; and for $i\geq 3$, $w_{-i}x^i=(h-1) \cdot (h+1)\cdots (h+i-2)\cdot (h+i)$ and $x^iw_{-i}=
(h-i-1) \cdot (h-i+1)\cdots (h-2)\cdot h$. 

The equalities in statement 2 are obvious. Then, by statement 1,  the elements $h$ and  $w_i$ ($i=\pm 1, \pm 2, \pm 3$) are algebra generators of $\CD (A)$. The subalgebra $\L :=K\langle h, w_3, w_{-3}\rangle $ is a generalized Weyl algebra $D[w_3, w_{-3}; \s , a=
(h+3)(h+1)(h-1)]$ which is a Noetherian algebra, \cite{GWA-AiA}. The algebra $\CD (A)$ is a finitely generated left and right $\L$-module, hence $\CD (A)$ is Noetherian.

3. By statement 2, $\Der_K(A)=K[x]h$ since $w_1=xh$. The rest follows.

4. The Jacobian ideal $\ga_1$ of $A$ is  $\D (R)$-stable since $\Der_K(A)= K[x]h$ and $h (\ga_1)\subseteq \ga_1$. Since $x^2\in \ga_1$ and $\o_{-2}(x^2)=-2\not\in \ga_1$, so the ideal $\ga_1$ is not $\CD (R)$-stable.  $\Box $

{\small

Department of Pure Mathematics

University of Sheffield

Hicks Building

Sheffield S3 7RH

UK

email: v.bavula@sheffield.ac.uk
}

\end{document}